\documentclass[12pt]{amsart}
\usepackage{amsmath,amsfonts,amssymb,amsthm}
\newtheorem{theorem}{Theorem}[section]
\newtheorem{lemma}[theorem]{Lemma}
\newtheorem*{lemmaA}{Lemma A.1}
\newtheorem{corollary}[theorem]{Corollary}

\newtheorem{definition}[theorem]{Definition}
\newtheorem{proposition}[theorem]{Proposition}
\newcommand{\R}{\mathbb{R}}             %% Real numbers
\newcommand{\C}{\mathbb{C}}             %% Complex numbers
\newcommand{\Sk}{{\mathcal S}}
\newcommand{\nuba}{\overline{\nu}}
\newcommand{\Om}{\Omega}
\newcommand{\Na}{\nabla}
\newcommand{\dinv}{\partial_{s}^{-1}}
\newcommand{\dnmap}{\Lambda_{\gamma}}
\newcommand{\bdom}{\del\Om}

\newcommand{\barg}{\overline{g}}
\newcommand{\Omc}{\Omega ^c}

\newcommand{\ol}{\overline} 
\newcommand{\del}{\partial}
\newcommand{\dba}{\ol\del}

\DeclareMathOperator{\re}{Re}

\title{Reconstruction of less regular conductivities in the plane}
\author{Kim Knudsen}
\address{\hskip-\parindent
Kim Knudsen\\
Department of Mathematical Sciences\\
Aalborg University\\
 Fredrik Bajers Vej 7G\\
 DK-9220 Aalborg \O\\
 Denmark.}
\email{kim@math.auc.dk}

\author{Alexandru Tamasan}
\address{\hskip-\parindent
Alexandru Tamasan\\
Department of Mathematics\\
 University of Washington\\
 Box 354350\\
 Seattle, Wa 98195-4350}
\email{tamasan@math.washington.edu}

\thanks{Research at MSRI is supported in part by Aalborg University,
Denmark, and in part by NSF grant DMS-0070488.}

\begin{document}
\begin{abstract}
We study the inverse conductivity problem of how to reconstruct an
isotropic electrical conductivity distribution $\gamma$ in an object
from static electrical measurements on the boundary of the object. We
give an exact reconstruction algorithm for the conductivity $\gamma\in
C^{1+\epsilon}(\ol \Om)$ in the plane domain $\Omega$ from the
associated Dirichlet to Neumann map on $\partial \Om.$ Hence we
improve earlier reconstruction results. The method used relies on a
well-known reduction to a first order system, for which the
$\ol\partial$-method of inverse scattering theory can be applied.
\end{abstract}
\maketitle

%\tableofcontents

\section{Introduction}
Let $\Omega\subset \R^n$ be a bounded and smooth domain and let
$\gamma \in L^\infty(\Omega)$ with $0 < c \leq \gamma(x), \text{
a.e. }x \in \Omega$ for some constant $c.$ Consider the conductivity
problem
\begin{align}\label{condprob}
\begin{split}
\nabla\cdot\gamma\nabla u = 0, \text{ in } \Omega,\\
u = f, \text{ on } \partial \Omega,
\end{split}
\end{align}
 for $f\in H^{1/2}(\partial \Omega).$ This problem has a unique
 solution $u\in H^1(\Omega)$ and hence we can define the Dirichlet to
 Neumann map
\begin{align*}
\Lambda_\gamma: H^{1/2}(\partial \Omega) &\rightarrow
 H^{-1/2}(\partial \Omega),\\  f &\mapsto \gamma \frac{\partial
 u}{\partial\nu}|_{\partial \Omega},
\end{align*}
where $\nu$ is the outer unit normal defined in a neighborhood
of $\partial \Omega.$ This definition should be understood in the
sense, that if $f,g \in H^{1/2}(\partial \Omega),$ then
\begin{align*}
  \langle\Lambda_\gamma f,g\rangle = \int_\Omega \gamma \nabla u \cdot
  \nabla v dx,
\end{align*}
where $ u \in H^1(\Omega) $ solves \eqref{condprob} and $v\in
H^1(\Omega)$ has $v|_{\partial \Omega}=g.$

%We will frequently use the fact, that when
%$\gamma\in C^{1+\epsilon}(\ol \Omega)$ then $\Lambda_\gamma :
%C^{1+\epsilon}(\partial \Omega) \rightarrow C^{\epsilon}(\partial
%\Omega).$

In this paper we study the problem of reconstructing $\gamma$ from
knowledge of $\Lambda_\gamma.$ This is the reconstruction issue of the
inverse conductivity problem.

The main motivation for studying this problem comes from a method for
medical imaging called Electrical Impedance Tomography (EIT) (see
\cite{cheney} for a thorough review of current methods in EIT). In EIT
one has a conductive body $\Omega$ with unknown conductivity $\gamma,$
and from static electric measurements on the boundary, i.e. by
applying a voltage potential on $\partial \Omega$ and measuring the
current flux through $\partial \Omega,$ one would like to monitor the
interior conductivity. Since tissue, bones, lungs etc. have different
conductive properties, an image of the conductivity is an image of the
internal structures.

The mathematical formulation of the inverse conductivity problem was
introduced by Calder\'on in \cite{calderon}. Calder\'on gave an
approximate reconstruction method for $\gamma$ close to
constant. Since Calder\'on's fundamental paper a lot of work has been
done on this problem.

The first global result was in \cite{kohnVogelius1}, where it was
proved, that if $\gamma\in C^\infty(\ol\Om)$ then $\Lambda_\gamma$
determines $\gamma|_{\partial \Omega}$ and all partial derivatives on
$\partial \Omega.$ From here it follows by analytic continuation, that
if $\gamma$ is real-analytic in $\Omega$ then it is determined by
$\Lambda_\gamma.$ Later in \cite{kohnVogelius2} the result was
extended to include piecewise real-analytic conductivities.

In three dimensions a major breakthrough in the theory of inverse
boundary value problems was the paper
\cite{sylvesterUhlmann87}. There it was proved that $\Lambda_\gamma$
determines $\gamma \in C^\infty(\ol \Omega).$ Their method of proof
relies on the existence of special solutions to the conductivity
equation, the so-called complex geometrical optics solutions or
exponentially growing solutions, which behave like $e^{x\cdot\zeta}$
for large $|x|$ and for $\zeta \in \C^n$ with $\zeta\cdot\zeta = 0.$
The assumption on $\gamma$ has since been relaxed by a number of
authors
(\cite{alessandrini88},\cite{nachman88},\cite{NSU},\cite{brown96}). The
best uniqueness result so far is for $\gamma \in C^{3/2}(\ol \Omega)$
obtained in \cite{PPU}.  Using exponentially growing solutions Nachman
(\cite{nachman88}) gave a reconstruction algorithm for conductivities
$\gamma \in C^{1,1}(\ol \Omega).$ Furthermore, Alessandrini
(\cite{alessandrini88}) proved a conditional stability result.

In two dimensions, uniqueness for small $\gamma \in
W^{3,\infty}(\Omega)$ was proved in \cite{sylvesterUhlmann86}, but the
general question remained open until 1996, when Nachman
\cite{nachman96} gave a uniqueness proof and a reconstruction
algorithm for $\gamma \in W^{2,p}(\Omega), p>1.$ Nachman's method
relies on the reduction of the conductivity equation to a
Schr\"odinger equation, which essentially requires two derivatives of
the conductivity.  For this equation he proved existence of
exponentially growing solution, which can be reconstructed from
$\Lambda_\gamma$ by the so-called $\ol \partial$-method in inverse
scattering theory originally developed by Beals and Coifman (see
\cite{bealscoifman2} and the reference there). Then using the
exponentially growing solutions Nachman was able to reconstruct
$\gamma.$ The algorithm has been tried out numerically (see
\cite{siltanenMuellerIsaacson}). In 1997 Brown and Uhlmann
(\cite{brownUhlmann}) improved the uniqueness result for $\gamma \in
W^{1,p}(\Omega), p>2.$ Their method relies on a reduction to a first
order system for which only one derivative is required. For this
system the $\ol\partial$-method is again applicable, see Section
\ref{outline} for further details. For both works conditional
stability has been proved, see \cite{liu} and \cite{barceloRuiz}.

For a thorough review on the history and developments of the inverse
conductivity problem we refer to \cite{uhlmann}.
 
In this paper we will show how the uniqueness proof in
\cite{brownUhlmann} for the two-dimensional inverse conductivity
problem can be turned into a reconstruction method. The reconstruction
algorithm parallels Nachman's algorithm, but since the reduction used
here requires less regularity, we can improve the reconstruction
result in \cite{nachman96}.  The main result is the following:
\begin{theorem}\label{recThm}
  Let $\Omega \subset \R^2$ be a bounded and smooth domain. Let
  $\epsilon>0$ and let $\gamma \in C^{1+\epsilon}(\ol \Omega)$ be such
  that $0<c\leq \gamma(x)$ for some constant $c.$ Then $\gamma$ can be
  reconstructed from $\Lambda_\gamma.$
\end{theorem}

The proof is given in the following sections as follows: In Section
\ref{outline} we introduce notation, outline the main idea of the
uniqueness proof and give the reconstruction algorithm. The steps of
the algorithm are then justified in the following sections.

\section{Outline of the reconstruction}\label{outline}
In this section we introduce notations and review the reduction of
\eqref{condprob} to a first order system. We then formulate an inverse
problem for this system and finally outline the reconstruction method.

Generally we work under the hypothesis stated in Theorem \ref{recThm}.
Furthermore, we assume that $\gamma = 1$ near $\del\Om.$ It was
noticed by Nachman in \cite{nachman96}, that this assumption is not a
restriction. When $\gamma \in C^{1+\epsilon}(\ol \Om)$ it is
well-known
(\cite{alessandrini88},\cite{sylvesterUhlmann88},\cite{nachman88})
that $\gamma$ and its normal derivative on the boundary are determined
by $\dnmap.$ Hence $\gamma$ can be extended outside $\Om$, such that
it is identically $1$ outside a compact set and $\gamma\in
C^{1+\epsilon}(\R^2)$. Moreover, the Dirichlet to Neumann map on the
boundary of a larger domain containing the original domain can be
reconstructed from the original Dirichlet to Neumann map and from the
known extension. We show in the appendix how to do such a construction
for $\gamma \in L^\infty (\Om)$.

%The following notations are taken from \cite{brownUhlmann}.  The
%matrix operator $D$ and the matrix potential $Q$ are defined by
%\begin{align}
%  D = 
%  \begin{pmatrix}
%    \dba & 0 \\
%    0 & \del
%  \end{pmatrix},
%\qquad
%  Q = 
%  \begin{pmatrix}
%    0 & q \\
%    \ol q & 0 
%  \end{pmatrix}, 
%\end{align}
%where $\del =(\del_x-i\del_y)/2$ and $\dba=(\del_x+i\del_y)/2$.
%For any $2\times 2$ matrix $A$, we let
%\begin{align}
%A^{off}=
%\begin{pmatrix}
%0 & a_{12}\\
%a_{21} & 0
%\end{pmatrix},
%\qquad
%A^{d}
%\begin{pmatrix}
%a_{11} & 0\\
%0 & a_{22}
%\end{pmatrix},
%\end{align}
%be the off-diagonal and diagonal parts of the matrix.

%For $z,k\in\C$ we consider the oscillatory maps 
%$e(z,k)=\exp(izk +i\ol{zk})$ and define the matrix
%\begin{align}
%  \Lambda_k(z) &= 
%  \begin{pmatrix} e(z, \ol k) &0 \\ 0 & e(z,-k)
%  \end{pmatrix}.
%\end{align}
%We let $E_k$ be the operator acting on $2\times 2$ matrices $A$ as
%\begin{align}
%E_k A &= A^{d}+\Lambda_k^{-1}A^{off}\\
%&= A^{d}+A^{off}\Lambda_{\ol k}
%\end{align}

It is well-known from the theory of elliptic equations, that if $f \in
C^{2+\epsilon}(\partial \Omega)$ and $\gamma \in C^{1+\epsilon}(\ol
\Omega), $ then \eqref{condprob} has a unique solution $u \in
C^{2+\epsilon}(\ol \Omega).$ Then the vector valued function $(v,w)^T$
defined by
\begin{align}\label{reduction}
  \begin{pmatrix} v \\ w \end{pmatrix} = \gamma^{1/2} \begin{pmatrix}
    \del u \\ \dba u \end{pmatrix} \in C^{1+\epsilon}(\ol\Omega)
    \times C^{1+\epsilon}(\ol\Omega)
\end{align}
satisfies the equation
\begin{align}\label{Dsystem}
  (D-Q)
  \begin{pmatrix}
    v \\
    w   
  \end{pmatrix} 
  = 0.
\end{align}
Here the matrix operator $D$ and the matrix potential $Q$ are defined
by
\begin{align}\label{DQ}
  D = 
  \begin{pmatrix}
    \dba & 0 \\
    0 & \del
  \end{pmatrix},
\qquad
  Q = 
  \begin{pmatrix}
    0 & q \\
    \ol q & 0 
  \end{pmatrix}, 
\end{align}
where $\del =(\del_x-i\del_y)/2$ and $\dba=(\del_x+i\del_y)/2$ and
\begin{align}\label{qgamma}
  q=-\gamma^{-1/2}\del\gamma^{1/2}.
\end{align}
%\begin{align}\label{Qcond} Q =
%-\gamma^{-1/2} \begin{pmatrix} 0 & \del\gamma^{1/2} \\
%\dba\gamma^{1/2}& 0 \end{pmatrix}.
%\end{align}

We shall see below in the proof of Theorem \ref{equivalence} that
there is a one-to-one correspondence between $u$ and $(v,w)^T,$ i.e.
given a solution $(v,w)$ to \eqref{Dsystem}, where $Q$ is of the form
\eqref{DQ} with $q$ given by \eqref{qgamma}, there exists a solution
$u$ to \eqref{condprob} such that \eqref{reduction} holds.

With $Q \in C^\epsilon(\R^2)$, the natural Cauchy data for the system
\eqref{Dsystem} is
\begin{align*}
{\mathcal C}_Q = \{(v|_{\bdom} , w|_{\bdom}): (v,w)\in
C^{1+\epsilon}(\ol\Om) \times C^{1+\epsilon}(\ol \Om) ,\, (D-
Q)(v,w)^T =0\}.
\end{align*}
We can now pose the following inverse problem for the system
\eqref{Dsystem}: Let $q\in C^{\epsilon}(\ol\Omega)$ and $Q$ defined by
\eqref{DQ}. Is $Q$ uniquely determined from $\mathcal C_Q$ and, in
case it is, how can it be reconstructed?

In this paper we show how to reconstruct $Q$ coming from a
conductivity equation. For this particular case we give an explicit
characterization of the Cauchy data ${\mathcal C}_Q$ in terms of the
Dirichlet to Neumann map $\Lambda_\gamma$ and this characterization
enables us to solve the inverse problem.

Brown and Uhlmann (\cite{brownUhlmann}) were the first to use
\eqref{Dsystem} in the study of the inverse conductivity problem and
their uniqueness proof was based on the scattering theory for this
system developed by Beals and Coifman (\cite{bealsCoifman}).  For
the direct scattering, one looks for a family of functions, which
solve the equation \eqref{Dsystem} in the whole plane (with $Q$
extended by zero outside $\Omega$) and have a certain exponential
behaviour at infinity.  More precisely, for each $k\in \C$, one seek
solutions $\Psi(\cdot, k)$ of the form
\begin{align}\label{psi}
\Psi(z,k) &= m(z,k) 
\begin{pmatrix}
e^{izk} & 0\\
0 & e^{-i\ol z k}
\end{pmatrix},
\end{align}
where $m$ approaches the identity matrix $I$ as $|z| \rightarrow
\infty$ and $z = x_1 + i x_2 \in \C.$ These solutions are called
exponentially growing solutions or complex geometrical optics.

A simple calculations shows that $m$ should solve the equation
\begin{align}\label{Dksystem}
  (D_k - Q) m = 0, 
\end{align} 
where $D_k$ is the matrix operator
\begin{align*}
  D_k = 
  \begin{pmatrix}
    (\ol\del-i\ol k) & 0 \\
    0 & (\del+ik) 
  \end{pmatrix}.
\end{align*}

The analysis of the direct scattering problem along these lines is due
to Beals and Coifman (\cite{bealscoifman3},\cite{bealsCoifman}). A
detailed analysis for the case $Q\in L^1(\R^2)\cap L^\infty(\R^2)$ can
be found in \cite{sung1},\cite{sung2},\cite{sung3}, and the analysis
for the case of compactly supported $Q\in L^p(\R^2),p>2,$ can be found
in \cite{brownUhlmann2}. Given our assumptions on $Q$ we have the
following theorem, see \cite{brownUhlmann2} and \cite{barceloRuiz}.
\begin{theorem}\label{direct}
Let $Q$ be a hermitian, compactly supported matrix with $Q \in
C^{\epsilon}(\R^2).$ Then for any $2<r<\infty$ and $k\in \C$ there is
a unique solution $m(z,k)$ to \eqref{Dksystem} with $(m(\cdot,k)-I)
\in L^{r}(\R^2)\cap C^{\epsilon}(\R^2).$
 
Furthermore, for fixed $z \in \R^2,$ the map $k\mapsto m(z ,k)$ is
differentiable, and it satisfies the equation
\begin{align}\label{dbar}
  \frac{\del}{\del \ol k}m(z,k) = m(z,\ol k)\Lambda_k(z)S(k),
\end{align}
where
\begin{equation}\label{S} 
  S(k)= \frac{i}{\pi}\int_{\C}
       \begin{pmatrix} 
          0 & q(z) m_{22}(z,k)\\ 
          -\ol q(z) m_{11}(z,k) & 0
        \end{pmatrix} 
      d\mu(z)
\end{equation}
and
\begin{align*}
  \Lambda_k(z) &= 
  \begin{pmatrix} e(z, \ol k) &0 \\ 0 & e(z,-k)
  \end{pmatrix}
\end{align*} 
for $e(z,k)=\exp(izk +i\ol{zk}).$
\end{theorem}
In the formula \eqref{S} $m_{11}$ and $m_{22}$ are the diagonal
components of $m$, and $d\mu(z)$ is the usual area measure. Note, that
since $Q \in C^\epsilon(\R^2),$ elliptic regularity theory gives, that
$m\in C^{1+\epsilon}_{loc}(\R^2).$ This fact will be used later.

The matrix $S(k)$ is called the scattering matrix.  Using the
definition of $\Psi$ and the fact that $Q$ is supported in $\Om$, an
integration by parts in \eqref{S} shows that
\begin{align}
  S(k)= \frac{i}{2\pi}\int_{\del\Om}
        \begin{pmatrix} 0 &   e^{-iz\ol{k}}\nu(z) \psi_{12}(z,k)\\
                          -e^{i\ol{ zk}}\ol \nu(z) \psi_{21}(z,k) & 0       
        \end{pmatrix} d\sigma(z), \label{Sform}
\end{align}
where $\psi_{12}$ and $\psi_{21}$ are the off diagonal entries of
$\Psi$.  In the above formula, $d\sigma(z)$ is the Lebesgue measure on
the boundary and $\nu$ is the complex normal at the boundary, i.e., if
$(\nu_1(z),\nu_2(z))$ denotes the outer unit normal at $z\in\bdom$,
then $\nu(z)= \nu_1(z)+i\nu_2(z)$ and $\nuba= \nu_1(z)-i\nu_2(z).$
Therefore, in order to determine the scattering matrix $S(k)$, it
suffices to find the traces on $\del\Om$ of the exponentially growing
solutions.

Once $S(k)$ is known, we can solve the $\ol \del$-equation
\eqref{dbar} to find $m(z,k)$ for each fixed $z\in\C$ and reconstruct
$Q$ by the formula
\begin{align*}
  Q(z) = \lim_{k_0 \rightarrow \infty}
  \mu(B_r(0))^{-1}\int_{\{k:|k-k_0|<r\}} D_km(z,k)d\mu(k)
\end{align*}
from \cite{brownUhlmann}. This would solve the inverse problem for
\eqref{Dsystem}. Furthermore, $\gamma$ can be found by solving
\eqref{qgamma}.

Since we are mainly interested in reconstructing the conductivity, and
this approach then requires taking the parameter to infinity and
solving a boundary value problem, we shall take a different route. As
observed in \cite{barceloRuiz} we can reconstruct $\gamma$ directly
from $m(z,0).$ In fact we will see, that
\begin{align}\label{mgamma}
  \gamma^{1/2}(z) =  \re (\tilde m_{11}(z,0)+\tilde m_{21}(z,0)),
\end{align}
where now $\tilde m$ is the unique solution given by Theorem
\ref{direct} to $D_k\tilde m = -Q^Tm.$ It was noticed by Beals and Coifman 
that the scattering matrix $\tilde S$ associated the potential 
$-Q^T$ relates to the scattering matrix $S$ for $Q$ by the formula
\begin{align}\label{Srel}
  \tilde S(k) = S(-\ol k)^T.
\end{align}
Therefore $\tilde m$ can be found from $S$ by 
solving \eqref{dbar} with $S$ replaced by $ S(-\ol k)^T.$ 

We are now ready to outline the reconstruction method:
\begin{enumerate}
  \item Determine $\Psi$ on $\del \Om$ (see Theorem
  \ref{mainrezult} below) 
  \item Define the scattering matrix $S$ by
  \eqref{Sform} and $\tilde S$ by \eqref{Srel} 
  \item Solve the $\ol \del$-equation \eqref{dbar} with $\tilde S$
  for $\tilde m(z,k), z \in \Om$
  \item Recover $\gamma$ on $\Om$ from \eqref{mgamma}.
\end{enumerate}
 The implementation of this method is subject to further study.
%%%%%%%%%%%%%%%%%%%%%%%%%%%%%%%%%%%%%%%%%%%%%%%%%%%%%%%%%%%%%%%%%%
%%%%%%%%%%%%%%%%%%%%%%%%%%%%%%%%%%%%%%%%%%%%%%%%%%%%%%%%%%%%%%%%%%%
\section{A boundary relation}
%%%%%%%%%%%%%%%%%%%%%%%%%%%%%%%%%%%%%%%%%%%%%%%%%%%%%%%%%%%%%%%%%%%

In this section we give an explicit characterization of the Cauchy data
for the first order system \eqref{Dsystem}, in case the potential $Q$
comes from a conductivity.

Assume $u \in C^{2+\epsilon}(\ol \Om)$ solves \eqref{condprob} for
some $f \in C^{2 + \epsilon}(\del \Om)$. Then $(v,w)=\gamma^{1/2}(\del
u, \dba u)$ solves \eqref{Dsystem} and one easily finds  the relation
\begin{align}\label{forwardBR}
    \left.  \begin{pmatrix} v \\w \end{pmatrix}\right | _{\del \Om} =
    \frac{1}{2} \begin{pmatrix} \ol \nu & -i\ol\nu \\ \nu & i\nu
    \end{pmatrix} \begin{pmatrix} \Lambda_\gamma(f)\\\del_s (f)
    \end{pmatrix}, \text{ on } \partial \Om,
\end{align}
where $\partial_s$ denotes the tangential derivative on $\partial
\Omega.$ Inverting \eqref{forwardBR} gives
\begin{align}\label{forwardBR2}
    \left. 
    \begin{pmatrix}
        \Lambda_\gamma(f)\\\del_s (f)
    \end{pmatrix}=
    \begin{pmatrix}
         \nu & \ol\nu \\ i\nu &  -i\ol\nu
    \end{pmatrix}
    \begin{pmatrix}
        v \\w
    \end{pmatrix}\right | _{\del \Om}. 
\end{align}
Let $\dinv :C^{1+\epsilon}(\bdom)\to C^{2+\epsilon}(\bdom)$ denote the
normalized inverse of $\del_s$ on the boundary defined by
\begin{align*}
\dinv f(s)=\int_0^s fds + c(f), 
\end{align*}
where the measure $ds$ is the arc length on $\bdom $ and the constant
$c(f)$ is chosen so that $\int_{\bdom}\dinv f ds=0$. Then from
\eqref{forwardBR2} we arrive at the following boundary relation for
$(v,w)|_{\partial \Omega}$
\begin{equation*}
i \dnmap \dinv (\nu v|_{\partial \Omega} -\nuba w|_{\partial \Omega} ) = 
 (\nu v|_{\partial \Omega} +\nuba w|_{\partial \Omega}).
\end{equation*}

This relation motivates the following definition
\begin{definition}
  Let $(h_1 , h_2)$ be a pair of functions in $C^{1+\epsilon}(\bdom)\times
  C^{1+\epsilon}(\bdom)$. We say that $(h_1 , h_2)$ satisfies the
  boundary relation if
\begin{equation}\label{bdrelation}
i \dnmap \dinv (\nu h_1 -\nuba h_2 ) = (\nu h_1 +\nuba h_2).
\end{equation}
\end{definition}

We just saw that solutions to \eqref{Dsystem}, defined through
\eqref{reduction}, satisfy the boundary relation, but as the following
theorem shows, \eqref{bdrelation} is in fact a complete
characterization of ${\mathcal C}_Q:$
\begin{theorem}\label{equivalence}
If $Q$ is given by \eqref{qgamma} then 
\begin{align}
{\mathcal C}_Q =
\{ (h_1,h_2)\in C^{1+\epsilon}(\del\Om )\times C^{1+\epsilon}(\del\Om ):\,
i\dnmap \dinv (\nu h_1 -\nuba h_2 ) = \nu h_1 +\nuba h_2 \}.
\end{align}
\end{theorem}
\begin{proof}
  First we show that any pair $(h_1, h_2)\in C^{1+\epsilon}(\del\Om )
  \times C^{1+\epsilon}(\del\Om )$ which satisfies the boundary
  relation is in ${\mathcal C}_Q$.  Let $u\in C^{2+\epsilon}(\ol{\Om})$ be
  the unique solution of the Dirichlet problem
\begin{align*}
\begin{cases}
\Na \cdot\gamma \Na u=0, & \text{ in } \Om,\\ u = i\del_s^{-1}(\nu h_1
-\nuba h_2), & \text{ in }\del\Om .
\end{cases}
\end{align*}
Define a solution $(v,w)^T$ to \eqref{Dsystem} by the relation
\eqref{reduction} with $u$ from above.  That $(v,w)|_{\partial \Omega}
= (h_1,h_2)$ follows from \eqref{forwardBR} by the use of
\eqref{bdrelation}.

Conversely, we start with a solution $(v,w)\in C^{1+\epsilon}(\ol\Om
)\times C^{1+\epsilon}(\ol\Om )$ to \eqref{Dsystem} and check the
boundary relation for $(v, w)|_{\partial \Omega}$ .  Since $v$ and $w$
solve the linear system \eqref{Dsystem} and $Q$ is of the form
\eqref{qgamma}, we have the compatibility relation
\begin{align*}
\dba(\gamma^{-1/2}v) = \del(\gamma^{-1/2}w).
\end{align*}
The Poincar\'e lemma ensures now the existence of such a $u$ that
\begin{align*}
  \gamma^{-1/2}
  \begin{pmatrix} v \\ w \end{pmatrix} =
    \begin{pmatrix} \del u \\ \dba u \end{pmatrix}.
\end{align*}
It is easy to check that $u$ is a solution to the conductivity
equation in the form $2\del u\dba \gamma +2\dba u \del\gamma +
4\gamma\del\dba u =0$.  Now relation \eqref{forwardBR} with
$f=u|_{\partial \Omega}$ shows that $(h_1,h_2)$ solves the boundary
relation.
\end{proof}

\section{From $\Lambda_\gamma$ to $S(k)$}
In this section we show how to reconstruct the trace on $\bdom$ of the
exponentially growing solutions $\Psi(z,k)$ defined by \eqref{psi}.
Notice that the first row has entries which are analytic outside
$\Om$, while the second row has entries which are anti-analytic
outside $\Om$.  Moreover they have a prescribed behavior at infinity
and their traces on $\del \Om$ satisfy the boundary relation
\eqref{bdrelation}.  We will prove, that these relations are
sufficient to reconstruct the trace of $\Psi$ uniquely.

Due to the symmetries
\begin{align}\label{msymmetry}
  m_{11}(z,k) = \ol{m_{22}(z,\ol k)}, 
\qquad m_{21}(z,k) = \ol{m_{12}(z,\ol k)},
\end{align}
which follows from the differential equations, the asymptotic for the
columns of $m$ and the uniqueness in Theorem \ref{direct}, it suffices
to reconstruct the first column $(\Psi_{11},\Psi_{21})^T$ of
$\Psi(z,k)$.
%For the sake of brevity, we relabel
%$(\Psi_{11},\Psi_{21})^T$ by $(v,w)^T$.

For every $(\zeta,z)\in \C$ with $\zeta \neq z$ we introduce
$$g_k(\zeta,z)=\frac{1}{\pi}\frac{e^{-ik(\zeta-z)}}{\zeta-z},$$
a
Green's kernel for $\dba$ which also takes into account exponential grow
at infinity. Using $g_k(\cdot,z)$ with $z\in\del\Om$, we define the
single layer potentials $\Sk_k$ and $\ol{\Sk}_k$ as boundary integral
operators by
\begin{align}
\Sk_k f(z) = \int_{\bdom}f(\zeta)g_k(\zeta,z)d\zeta,\,\,\,
\ol{\Sk}_k f(z)= \int_{\bdom}f(\zeta)\barg_k(\zeta,z)d\zeta.
\end{align}
These integrals are understood in the sense of principal value. It is a 
classical result in singular integral theory that these operators are well 
defined, see \cite{musk}.
\begin{proposition}If $0<\epsilon< 1$ then 
  $\Sk_k ,\,\, \ol{\Sk}_k: C^{1+\epsilon}(\bdom) \to
  C^{1+\epsilon}(\bdom)$ are well defined bounded operators.
\end{proposition}

Let ${\mathcal K}$ denote the matrix operator
\begin{align*}
{\mathcal K}=\frac{1}{2}
   \begin{pmatrix}
    I-i\Sk _k&0\\
    0&I-\ol{\Sk}_k.
\end{pmatrix}
\end{align*}
The following lemma gives a necessary condition for a pair of
functions defined in $\Omc=\C\setminus \ol \Om$ to be analytic and
anti-analytic respectively and have a special exponential decay at
infinity:

\begin{proposition}
Let $v$ be analytic and $w$ be anti-analytic in $\Omc$ and continuous
to $\del\Om$, such that $ve^{-izk}-1\in L^r(\Omc)$ and $we^{-izk}\in
L^r(\Omc)$ for some $r\geq 1$. Then the trace $(h_1,h_2) =
(v,w)|_{\del \Om}$ solves
\begin{align}\label{infinity}
{\mathcal K}
\begin{pmatrix}
h_1\\
h_2
\end{pmatrix}
=\begin{pmatrix}
e^{izk}\\
0
\end{pmatrix}.
\end{align}
\end{proposition}
\begin{proof}
  We will prove the relation for $v;$ for $w$ the result follows by
  similar arguments. Let $z\in\Omc$ be a fixed point inside some large
  ball $B_R$ of radius $R$.  With the $\dba$ operator acting with
  respect to $\zeta$ and $R>0$ such that $\ol \Om \subset B(0,R),$
  Pompeiu's formula (\cite{vekua}) states, that 
\begin{align*}
  v(z)e^{-ikz}-1 &= \frac{1}{2\pi i}\int_{\del
    B_R}\frac{v(\zeta)e^{-ik\zeta}-1}{\zeta-z}d\zeta-
  \frac{1}{2\pi i}\int_{\del\Om}\frac{v(\zeta)e^{-ik\zeta}-1}{\zeta-z}d\zeta
  \\ &\qquad -\int_{B_R\setminus \ol \Omega}
  \frac{\dba(v(\zeta)e^{-ik\zeta}-1)}{\zeta-z}d\zeta  \\
  &= \frac{1}{2\pi i}\int_{\del
    B_R}\frac{v(\zeta)e^{-ik\zeta}-1}{\zeta-z}d\zeta-
  \frac{1}{2\pi i}\int_{\del\Om}\frac{v(\zeta)e^{-ik\zeta}-1}{\zeta-z}d\zeta.
\end{align*}
%But
%\begin{align*}
%\frac{1}{2i}\int_{\del B_R}v(\zeta)g_k(\zeta,z)d\zeta=
%\frac{e^{ikz}}{2\pi i}
%\int_{\del B_R}\frac{v(\zeta)e^{-ik\zeta}-1}{\zeta-z}d\zeta+
%e^{ikz}\frac{1}{2\pi i}\int_{\del B_R}\frac{1}{\zeta-z}d\zeta=\\
%\int_{\del B_R}\frac{v(\zeta)e^{-ik\zeta}-1}{\zeta-z}d\zeta+
%e^{ikz}.
%\end{align*}
%The last equality took into account the winding number of $\del B_R$ 
%around $z$. 
Since $ve^{-izk}-1 \in L^r(\R^2)$ the integral on $\del B_R$ converges
to zero at least on a sequence of increasing radii, and since $z\in
\Omc,$ Cauchy's formula gives
\begin{align*}
  \frac{1}{2\pi i}\int_{\del\Om}\frac{1}{\zeta-z}d\zeta = 0.
\end{align*}
Therefore
\begin{align}\label{vextension}
v(z)+\frac{1}{2i}\int_{\del\Om}v(\zeta)g_k(\zeta,z)d\zeta = e^{izk},
\,\, z\in\Omc.
\end{align}
Now let $z \in \Omc$ approach some point on the boundary and apply
Plemelj's Formula \cite{musk} to get the equation for $v$.
\end{proof}

We know that the exponentially growing solutions
$(\psi_{11}(\cdot,k),\psi_{21}(\cdot,k))|_{\partial \Omega}$ must
satisfy \eqref{infinity} and the boundary relation
\eqref{bdrelation}. The following result shows that these conditions
are sufficient to characterize the trace of the exponentially growing
solutions.
\begin{theorem}\label{mainrezult}
  The only pair $(h_1 ,h_2)^T\in C^{1+\epsilon}(\del\Om)\times
  C^{1+\epsilon}(\del\Om)$ which satisfies
\begin{align}
\begin{split} \label{defineK}
{\mathcal K}\begin{pmatrix}h_1\\h_2\end{pmatrix}
&=\begin{pmatrix}
e^{izk}\\0\end{pmatrix}\\
(I-i\dnmap\dinv)(\nu h_1)(z)&=(I+i\dnmap\dinv)(\nuba h_2)(z),
\end{split}
\end{align}
is $(\psi_{11}(\cdot,k),\psi_{21}(\cdot,k))|_{\partial \Omega}.$
\end{theorem}
\begin{proof}
  Let $h = (h_1,h_2) \in C^{1+\epsilon}(\del\Om)\times
  C^{1+\epsilon}(\del\Om)$ be any solution to \eqref{defineK}. Guided
  by \eqref{vextension} and the related formula for $w$ we extend $h$
  outside $\Om$ to $(v(z),w(z))$ by \begin{align*}
  v(z)=-\frac{1}{2i}\int_{\del\Om}h_1(\zeta)g_k(\zeta,z)d\zeta+
  e^{izk},\,\,\,
  w(z)=-\frac{1}{2i}\int_{\del\Om}h_2(\zeta)\barg_k(\zeta,z)d\zeta.
  \end{align*} Inside $\Om$ we extend $h$ as a solution to
  \eqref{Dsystem}, which can be done by Theorem \ref{mainrezult},
  since $h$ satisfies the boundary relation.
  
  We will prove, that $(v,w)$ is a solution to \eqref{Dsystem} and
  that $ve^{-izk}-1,w \in L^r(\R^2), \;r>2.$ Then the result follows from
  Theorem \ref{direct}.
  
  Outside $\Om$, $v$ is analytic and $w$ is anti-analytic. Moreover,
  since $\gamma$ is identically $1$ near the boundary, we also have
  that $v$ is analytic and $w$ is anti-analytic in an inside
  neighborhood of the boundary.  We need to check the analyticity of
  $v$ and anti-analyticity of $w$ as we go across the
  boundary. For this it suffices to have continuity across the
  boundary and then invoke Morera's theorem.  We shall show the
  continuity of $v$, for $w$ similar reasoning works.  Let $z$
  approach some point $z_0\in\del\Om$ from outside.  Then using
  Plemelj's formula we get
  \begin{align*}
    \lim_{z\to z_0^+}v(z)= -\left(-\frac{h_1(z_0)}{2}+ \frac{1}{2i}S_k
      h_1(z_0)\right)+e^{iz_0k}.
  \end{align*} 
  Now use the first equation of \eqref{defineK} to conclude
  $\lim_{z\to z_0}v(z)=h_1(z_0)$.  The continuity of $v$ from inside
  comes from the regularity theory for elliptic equations. 
  
  We are left to show that $ve^{-izk}-1 \in L^r(\R^2),$ but since
  $$v(z)e^{-izk}-1 = -\frac{1}{2\pi i} \int_{\del\Om}\frac {
h_1(\zeta)e^{-ik\zeta}}{\zeta -z}d\zeta= O\left( \frac{1}{|z|}\right),
\text{ as } |z|\to\infty,$$ it is in $L^r(\R^2)$ for $r>2$. This
proves the theorem.
\end{proof}

As an corollary we easily get
\begin{corollary}
$\ker {\mathcal K} \cap {\mathcal C}_Q=\{ 0\}$. 
\end{corollary}
\begin{proof} Let $(h_1, h_2)^T\in \ker {\mathcal K}\cap {\mathcal C}_Q$ and $(v,w)^T$ be the 
  trace of the exponentially growing solutions. Then $(h_1 + v, h_2 +
  w)^T$ is also satisfies the boundary relation and satisfies
  \eqref{defineK}. Now use the uniqueness of Theorem \ref{mainrezult}
  to conclude that $(h_1, h_2)^T =0$.
\end{proof}

To obtain the exponentially growing solutions on $\partial \Om$ one
has to solve the three equations \eqref{defineK}. Since
$(e^{ikz},0)^T$ solves
\begin{align*}
{\mathcal K}\begin{pmatrix} e^{izk}\\0\end{pmatrix} &=\begin{pmatrix}
e^{izk}\\0\end{pmatrix}, 
\end{align*}
the pair of traces of the exponentially growing solutions can be
constructed as the projection of $(e^{ikz},0)^T$ on ${\mathcal C}_Q$ along
$\ker \mathcal K.$

\section{From $S(k)$ to $\gamma$}
In this section we will prove how to reconstruct $\gamma$ from $S(k).$
More precisely we will work with the dual $\tilde S(k) = S(-\ol k)^T,$
solve the equation \eqref{dbar} for $\tilde m$ and finally find
$\gamma$ by the formula \eqref{mgamma}.

Since the equation \eqref{dbar} is mixing terms, we will start by
writing pseudoanalytic equation for a suitable linear combination of
the elements in $m:$
\begin{lemma}
Let $z \in \C$ be fixed and let 
\begin{align}\label{mp}
  \tilde m_+(z,k) = \tilde m_{11}(z,k) + \ol {\tilde m_{12}(z, \ol
    k)}.
\end{align} Then for $k \in \C$ we
have the equation
\begin{align}
  \frac{\partial}{\partial \ol k} \tilde m_+ (z,k) &=
  e(z,-k)S_{21}(k)\ol{\tilde m_+(z, k)}\label{pseudo1},% \\
%  \frac{\partial}{\partial \ol k} \tilde m_- (z,k) &=
%  -e(z,-k)S_{21}(k)\ol{\tilde m_-(z, k)}. \label{pseudo2}
\end{align}
\end{lemma}
\begin{proof}
  The result follows from \eqref{dbar} and \eqref{msymmetry}.
\end{proof}

To solve this equation we will need a result about the behavior of
$\tilde m_+$ for large $k:$
\begin{theorem}[\cite{brownUhlmann}]\label{kasymp}
  Suppose $Q\in L^p(\R^2), p>2$ is compactly supported. Then for all $q>2p/(p-2)$
  \begin{align*}
    \sup_{z \in \C} \|m(z,\cdot)-I\|_{L^q(\R^2)} \leq C
  \end{align*}
  where $C$ depends on $p,q$ and $Q.$ 
\end{theorem}
Furthermore, we will need a certain decay of $S.$ In
\cite{brownUhlmann} the fact that $S \in L^2(\R^2)$ was sufficient to
prove uniqueness of solutions to \eqref{pseudo1} subject to the
asymptotic given in Theorem \ref{kasymp}. We will however, need a
slightly better decay for $S$ in order to write the integral equations
corresponding to \eqref{pseudo1}. The following lemma from
\cite{barceloRuiz} gives the result:
\begin{lemma}\label{Sreg}
  Let $Q \in C^{\epsilon}(\R^2)$ be compactly supported. Then $S \in
  L^r(\R^2)$ for any $r>r_0= 4/(2+\epsilon).$
\end{lemma}

To simplify notations we introduce the operator $(\partial/\partial
\ol k)^{-1} = \partial_{\ol k}^{-1}$ defined by
\begin{align*}
  \partial_{\ol k}^{-1} f (k) = \frac{1}{\pi} \int_\C \frac{f(k')}{k-k'}dk'.
\end{align*}
 We can now write the integral equations for $\tilde m_+$ and prove
unique solvability of these:
\begin{lemma}
  Let $Q\in C^{\epsilon}(\R^2)$ be compactly supported, and let $z \in
  \C$ be fixed. Then for $q>4/\epsilon$ the equations 
\begin{align}
  (I - \partial_{\ol k}^{-1}(e(z,-k)S_{21}\cdot))(\tilde m_+-I) &=
  \partial_{\ol k}^{-1}(e(z,-k)S_{21}), \label{intm1} %\\ (I +
%  \partial_{\ol k}^{-1}(e(z,- k)S_{21}\cdot))(\tilde m_--1) &=
%  \partial_{\ol k}^{-1}(e(z, -k)S_{21}). \label{intm2}
\end{align}
have the unique solution $\tilde m_+(z,\cdot)-I$ in $L^q(\R^2).$
Furthermore, $\tilde m_+(z,\cdot)-I \in C^\alpha(\R^2),$ for $\alpha <
(1+\epsilon)/2.$
\end{lemma}
\begin{proof}
  Since $S_{21}\in L^2(\R^2)$ and $|e(z,-k)|=1$ we know (see
  \cite[Lemma 4.2]{nachman93}) that $\partial_{\ol
  k}^{-1}(e(z,-k)S_{21}\cdot)$ is a compact operator in $L^s(\R^2)$
  for $2<s<\infty.$ Furthermore, since $S_{21} \in L^r(\R^2)$ for some
  $4/(2+\epsilon)<r<2$ (Lemma \ref{Sreg}), it follows by the
  Hardy-Littlewood-Sobolev inequality (\cite{stein}) that
  $\partial_{\ol k}^{-1}(e(z,-k)S_{21}) \in L^q(\R^2)$ for
  $q>\epsilon/2.$ Hence the integral equation \eqref{intm1} in
  $L^q(\R^2)$ follows from \eqref{pseudo1} and Lemma \ref{kasymp}.
  
  The unique solvability of \eqref{intm1} in $L^q(\R^2)$ follows from
  Fredholm alternative if $(I + \partial_{\ol
    k}^{-1}(e(z,-k)S_{21}\cdot))$ has trivial kernel in $L^q(\R^2),$
  but this is known from Liouville's theorem for pseudoanalytic
  functions since $e(z,-k)S_{21} \in L^r(\R^2)\cap L^{r'}(\R^2)$ (see
  \cite{brownUhlmann},\cite{sung1},\cite{nachman96},\cite{vekua}).

  To prove the H\"older continuity of $\tilde m_+(z,\cdot)-I$ we use the
  fact, that convolution by $1/z$ maps $L^p(\R^2)$ into
  $C^\alpha(\R^2)$ for $1<p<\infty,\; \alpha=1-1/{p_0}$ and $
  \max(2,p)< p_0<\infty,$ see \cite{sunUhlmann}.
\end{proof}

By the previous theorem we find for fixed $z \in \C$
\begin{align*}
  (\tilde m_+(z,k)-I) &= (I - \partial_{\ol
    k}^{-1}(e(z,-k)S_{21}\cdot))^{-1}(\partial_{\ol
    k}^{-1}(e(z,-k)S_{21})),% \\ (\tilde m_-(z,k)-I) &= (I + \partial_{\ol
%    k}^{-1}(e(z,- k)S_{21}\cdot))^{-1}(\partial_{\ol
%    k}^{-1}(e(z,-k)S_{21})),
\end{align*}
in $L^q(\R^2) \cap C^\alpha(\R^2)$ for $\alpha = 1-(1/r+1/q).$ To
reconstruct $\gamma$ we use the formula
\begin{align*}
  \gamma(z) &= \tilde m_{11}(z,0)+\ol{\tilde m_{21}(z,0)}
\end{align*}
from \cite[Proposition 4.2]{barceloRuiz}. Since $\gamma$ is real, we
 rewrite the formula using \eqref{msymmetry} and \eqref{mp}
\begin{align*}
  \gamma(z) &= \tilde m_{11}(z,0)+\ol{\tilde m_{21}(z,0)}\\ 
  &= \re\tilde m_{11}(z,0) + \re \ol{\tilde m_{21}(z,0)}\\
  &= \re\tilde m_{11}(z,0) + \re {\tilde m_{12}(z,0)}\\
  &=\re \tilde m_+(z,0).
\end{align*}
This ends the reconstruction.

%%%%%%%%%%%%%%%%%%%%%%%%%%%%%%%%%%%%%%%%%%%
%%%%%%%%%%%%%%%%%%%%%%%%%%%%%%%%%%%%%%%%%%%%%%
%%APPENDIX
%%
%%%%%%%%%%%%%%%%%%%%%%%%%%%%%%%%%%%%
%%%%%%%%%%%%%%%%%%%%%%%%%%%%%%%%%%%%

\section*{Appendix}
We will prove, that the Dirichlet to Neumann map for an extended
domain $\Omega'$ can be calculated from the Dirichlet to Neumann map
for the domain $\Omega$ given the value of the conductivity in
$\Omega'\setminus \ol \Omega:$
\begin{lemmaA}
Let $\gamma \in L^\infty(\Om)$ and let $\Lambda_\gamma$ be
given. Let $\Om'\subset \R^2$ be such that $\ol \Om \subset
\Om'.$ Take any function $\gamma' \in L^\infty(\Om')$ with
$\gamma'|_{\Om} = \gamma.$ 

Then $\Lambda_{\gamma'}$ can be recovered from $\Lambda_\gamma$ and
$\gamma'|_{\Om'\setminus{\ol \Om}}$
\end{lemmaA}
\begin{proof}
  From the definition of the Dirichlet to Neumann map we have for any
  $f,g \in H^{1/2}(\del \Om'),$ that
\begin{align*}
  \langle\Lambda_{\gamma'}f,g\rangle &= \int_{\Om'}\gamma'\nabla u \cdot \nabla v\\
  &= \int_{\Om'\setminus\ol\Om} \gamma'\nabla u \cdot \nabla v +   
  \langle\Lambda_{\gamma}u|_{\del \Om},v|_{\del \Om}\rangle,
\end{align*}
where $u \in H^1(\Om')$ denotes the unique solution to
\begin{align}\label{extProb}
\begin{split}
\nabla\cdot \gamma'\nabla u &= 0, \text{ for } x \in \Om',\\
u &= f, \text{ for } x \in \del \Om',
\end{split}
\end{align}
and $v \in H^1(\Om')$ is any function with $v|_{\del \Om'} =
g.$ Hence we see, that $\Lambda_{\gamma'}$ can be calculated from
$\gamma|_{\Om'\setminus\ol\Om}$ and $\Lambda_\gamma$ without
explicit knowledge of $\gamma$ in $\Om,$ provided the solution $u$
to \eqref{extProb} can be found in $\Om'\setminus \Om.$

We claim, that $u$ in $\Om'\setminus \Om$ can be found as the unique solution to 
\begin{align}\label{redProb}
\begin{split}
\nabla\cdot \gamma'\nabla u &= 0, \text{ for } x \in \Om'\setminus
\ol\Om,\\ u &= f, \text{ for } x \in \del \Om',\\
\gamma'\frac{\del u}{\del \nu}|_{\del
\Om}&=\Lambda_\gamma(u|_{\del \Om}), \text{ for } x \in
\del \Om'.
\end{split}
\end{align}

That $u|_{\Om'\setminus\ol\Om}$ solves \eqref{redProb} follows,
since for any $v \in H_0^1(\Om')$
\begin{align*}
  \int_{\Om'\setminus\ol\Om} \gamma'\nabla u \cdot \nabla v +
  \langle\Lambda_{\gamma}u|_{\del \Om},v|_{\del
  \Om}\rangle &= \int_{\Om'} \gamma' \nabla u \cdot \nabla v \\
  &=0.
\end{align*}
For the uniqueness we assume, that $u_0\in H^1(\Om'\setminus \ol
\Om)$ solves \eqref{redProb} with $f=0.$ Now extend $u_0$ into
$\Om$ as the solution to
\begin{align*}
\begin{split}
\nabla\cdot \gamma\nabla u &= 0, \text{ for } x \in \Om,\\
u &= u_0, \text{ for } x \in \del \Om.
\end{split}
\end{align*}
But since $u_0$ is in $H^1(\Om')$ and solves \eqref{extProb} with
$f=0,$ we conclude, that $u_0=0$ in $\Om'.$
\end{proof}

%%%%%%%%%%%%%%%%%%%%%%%%%%%%%%%%%%%%%%%%%%%%%%%%%%%%%%%%%%%%%%%%%%%%%%

%\bibliographystyle{amsalpha}

%\bibliography{invCond}

\end{document}